\documentclass{amsart}
\newtheorem{thm}{Theorem}[section]
\newtheorem{cor}[thm]{Corollary}
\newtheorem{conj}[thm]{Conjecture}

\newtheorem{lem}[thm]{Lemma}
\newtheorem{prop}[thm]{Proposition}
\theoremstyle{definition}

\theoremstyle{remark}
\newtheorem{rem}[thm]{Remark}
\numberwithin{equation}{section}

\begin{document}

\title[Advances on the Bessis-Moussa-Villani Trace Conjecture]
{Advances on the Bessis- \\ Moussa-Villani Trace Conjecture}%
\author{Christopher J. Hillar}%
\address{Department of Mathematics, Texas A\&M University, College Station, TX 77843}
\email{chillar@math.tamu.edu}
 \thanks{Supported under a National Science Foundation Postdoctoral Research Fellowship.}
  

\begin{abstract}
A long-standing conjecture asserts that the polynomial \[p(t) = \text{Tr}[(A+tB)^m]\] 
has nonnegative coefficients whenever $m$ is a positive integer and $A$ and 
$B$ are any two $n \times n$ positive semidefinite Hermitian matrices.
The conjecture arises from a question raised by Bessis, Moussa, and 
Villani (1975) in connection with a problem in theoretical physics.  
Their conjecture, as shown recently by Lieb and Seiringer, is equivalent
to the trace positivity statement above.  
In this paper, we derive a fundamental 
set of equations satisfied by $A$ and $B$ that minimize or maximize a 
coefficient of $p(t)$.
Applied to the Bessis-Moussa-Villani (BMV) conjecture,
these equations provide several reductions.  In particular, we prove that 
it is enough to show that (1) it is true for infinitely many $m$,
(2) a nonzero (matrix) coefficient of $(A+tB)^m$ always has at least one positive eigenvalue,
or (3)  the result holds for singular positive semidefinite matrices.  Moreover, we prove that if the conjecture is false for 
some $m$, then it is false for all larger $m$.
\end{abstract}
\maketitle

\section{Introduction}
In \cite{bessis}, while studying partition functions of quantum
mechanical systems, a conjecture was made regarding a positivity
property of traces of matrices.  If this property holds, explicit
error bounds in a sequence of Pad\'e approximants follow.
Let $A$ and $B$ be $n \times n$ Hermitian matrices with
$B$ positive semidefinite, and let 
\[ \phi^{A,B}(t) = \text{Tr}[\exp{(A-tB)}].\]  
The original formulation of the Bessis-Moussa-Villani conjecture asserts that
the function $\phi^{A,B}$ is completely monotone;
in other words, $\phi^{A,B}$ is the Laplace transform
of a positive measure $\mu^{A,B}$ supported on $[0,\infty)$:
\[ \text{Tr}[\exp{(A-tB)}] = \int_{0}^{\infty}\exp(-tx)d\mu^{A,B}(x).\]
Equivalently, the derivatives of the function $f(t) = \phi^{A,B}(t)$ alternate 
signs: \[ (-1)^n f^{(n)} (t) \geq 0, \ \ \ t > 0, \ n = 0,1,2,\ldots.\]

Since its introduction in \cite{bessis}, many partial results and 
substantial computational experimentation have been given 
\cite{KJC, KJC2,otherBMVhyp,otherBMVtry2, FH, JH, JH3, otherBMVtry1, Moussa}, all
in favor of the conjecture's validity.  
However, despite much work, very little is known about the
problem, and it has remained unresolved except in
very special cases.  Recently, Lieb and Seiringer in \cite{Lieb}, and 
as previously communicated to us
\cite{JH}, have reformulated the conjecture of \cite{bessis} as a
question about the traces of certain sums of words in two positive
definite matrices.  In what follows, we shall use the standard convention that
a \textit{positive definite} (resp. \textit{positive semidefinite}) matrix
is one that is complex Hermitian and has positive eigenvalues  
(resp. nonnegative eigenvalues).  
\begin{conj}[Bessis-Moussa-Villani]\label{BMV}
The polynomial $p(t) = \text{\rm{Tr}}\left[(A+tB)^m\right]$ has all
nonnegative coefficients whenever $A$ and $B$ are $n \times n$ positive
semidefinite matrices.
\end{conj}
\begin{rem}
Although not immediately obvious, the polynomial $p(t)$
has all real coefficients (see Corollary \ref{realcoeffcor}).
\end{rem}

The coefficient of $t^k$ in $p(t)$ is the trace of $S_{m,k}(A,B)$,
the sum of all words of length $m$ in $A$ and $B$, in which $k$
$B$'s appear (it has been called the $k$-th Hurwitz product of $A$
and $B$).  In \cite{JH}, among other things, it was noted that,
for $m <6$, each constituent word in $S_{m,k}(A,B)$ has nonnegative
trace.  Thus, the above conjecture is valid for $m < 6$ and
arbitrary positive integers $n$.  It was also noted in \cite{JH} (see also \cite{bessis})
that the conjecture is valid for arbitrary $m$ and $n < 3$.  Thus,
the first case in which prior methods did not apply and the
conjecture was in doubt, is $m = 6$ and $n = 3$.  Even in this
case, all coefficients, except
$\text{Tr}[S_{6,3}(A,B)]$, were known to be nonnegative (also as shown
in \cite{JH}). It was only recently \cite{JH3}, using heavy computation,
that this remaining coefficient was shown to be nonnegative.  

Much of the subtlety of Conjecture \ref{BMV} can be seen by the fact 
that $S_{m,k}(A,B)$ need not have 
all nonnegative eigenvalues, and in addition that some words within the $S_{m,k}(A,B)$ 
expression can have negative trace (see \cite{JH}, where it is shown
that Tr$[ABABBA]$ can be negative).

Our advancement is the introduction of a fundamental pair of matrix 
equations satisfied by $A$ and $B$ that minimize or maximize
a coefficient of $p(t)$.  In what follows, we will be using the natural 
Euclidean norm on the set of complex $n \times n$ matrices: 
\begin{equation*}\label{euclidnorm}
 \|A\| = \text{Tr}[AA^*]^{1/2}.
\end{equation*} 
(Here, $C^*$ denotes the \textit{conjugate transpose} of a 
complex matrix $C$).  The precise statement of our main result is the following.
  
\begin{thm}\label{mainvarthmAB}
Let $m > k > 0$ be positive integers, and 
let $A$ and $B$ be positive semidefinite
matrices of norm $1$ that minimize (resp. maximize) the quantity 
$\text{\rm Tr}[S_{m,k}(A,B)]$
over all positive semidefinite matrices of norm $1$.
Then, $A$ and $B$ satisfy the following pair of equations:
\begin{equation}\label{matrixeqs}
\begin{cases}
AS_{m-1,k}(A,B) & = \ A^2 \text{\rm Tr}[AS_{m-1,k}(A,B)]  \\
BS_{m-1,k-1}(A,B) & = \ B^2 \text{\rm Tr}[BS_{m-1,k-1}(A,B)].
\end{cases}
\end{equation}
\end{thm}

We call (\ref{matrixeqs}) the \emph{Euler-Lagrange equations}
for Conjecture \ref{BMV}.  The name comes from the resemblance of 
our techniques to those of computing the first variation in the calculus 
of variations.  We should remark that there have been other variational 
approaches to this problem \cite{KJC,KJC2}; a review can be found in \cite{Moussa}.
Although we are motivated by Conjecture \ref{BMV},
we discovered that these equations are also
satisfied by a minimization (resp. maximization) over Hermitian matrices $A$ and $B$
of norm $1$ (see Corollary \ref{hermitianeqs}), and it is natural
to consider this more general situation.   In this regard, we present
the following application of the Euler-Lagrange equations.

\begin{thm}\label{maxtraceconj}
If $A$ and $B$ are Hermitian
matrices of norm $1$ and $m > 1$, then
\begin{equation*}
\left|\text{\rm Tr}[S_{m,k}(A,B)]\right| \leq {m \choose k}.
\end{equation*}
Moreover, if $m > k > 0$, then equality holds only when $A = \pm B$, and 
if in addition $m > 2$, then $A$ has precisely one nonzero eigenvalue.
\end{thm}
\begin{rem}
When $m = 1$, this theorem fails to hold.  For example, let 
$A$ be the $n \times n$ diagonal matrix $A = \text{diag}(n^{-1/2},\ldots,n^{-1/2})$.
Then $\|A\| = 1$, but $\text{\rm Tr}[S_{1,0}(A,B)] = \text{\rm Tr}[A] =  n^{1/2} > 1$ for 
$n > 1$.
\end{rem}

It is easy to see that this maximum is at least 
${m \choose k}$, and
using elementary considerations involving the Cauchy-Schwartz inequality,
one can show that
\[ \left|\text{Tr}[S_{m,k}(A,B)]\right| \leq \|S_{m,k}(A,B)\|n^{1/2} \leq 
{m \choose k}\|A\|^{m-k}\|B\|^k n^{1/2} = {m \choose k}n^{1/2}.\] 
However, we do not know if a dependency on the size of the matrices involved can
be removed without appealing to equations (\ref{matrixeqs}).   

As a strategy to prove Conjecture \ref{BMV},
we offer the following.

\begin{conj}\label{commuteconj}
Let $m > k > 0$ be positive integers.
Positive semidefinite (resp. Hermitian) matrices
$A$ and $B$ of norm $1$ that satisfy the Euler-Lagrange 
equations commute.
\end{conj}

From this result, Conjecture \ref{BMV} would be immediate.
Of course,  Theorem 
\ref{maxtraceconj} implies that Conjecture \ref{commuteconj} holds for 
the case of Hermitian matrices.
We next list some of the major consequences of the
equations found in Theorem \ref{mainvarthmAB}.  
The first one implies that counterexamples to Conjecture \ref{BMV} are
closed upwards.  The precise statement
is given by the following.

\begin{thm}\label{mainthm}
Suppose that there exist integers $M, K$ and $n \times n$ 
positive definite matrices $A$ 
and $B$ such that $\text{\rm Tr}[S_{M,K}(A,B)] < 0.$  Then, 
for any $m \geq M$ and $k \geq K$ such 
that $m-k \geq M-K$, there exist $n \times n$ positive definite
$A$ and $B$ making $\text{\rm Tr}[S_{m,k}(A,B)]$ negative.
\end{thm}

\begin{cor}
If the Bessis-Moussa-Villani conjecture is  true for some $m_0$, then it 
is also true for all $m < m_0$.
\end{cor}

This reduces the BMV conjecture to its ``asymptotic" formulation.

\begin{cor}
If the Bessis-Moussa-Villani conjecture is true for infinitely many $m$, then it 
is true for all $m$.
\end{cor}

A next result characterizes
the BMV conjecture in terms of the eigenvalues of the matrix $S_{m,k}(A,B)$.

\begin{thm}\label{poseigvaluethm}
Fix positive integers $m > k$ and $n$. Then,
\text{\rm Tr}$[S_{m,k}(A,B)] \geq 0$ for all $n \times n$ positive semidefinite $A$ and $B$  
if and only if whenever $S_{m,k}(A,B) \neq 0$, it has at least 
one positive eigenvalue.
\end{thm}

\begin{rem}
This theorem can be viewed as a transfer principle for the BMV conjecture:
instead of proving positivity for the sum of all the eigenvalues, we need only 
show it for at least one of them.
\end{rem}

Our final result generalizes 
a fact first discovered in \cite{JH3} (there only the real case was 
considered), and it implies that it is enough to prove the 
Bessis-Moussa-Villani conjecture for singular $A$ and $B$.

\begin{thm}\label{singthm}
Let $m,n$ be positive integers, 
and suppose that $\text{\rm{Tr}}\left[(A+tB)^{m-1}\right]$ has nonnegative
coefficients for each pair of $n \times n$
positive semidefinite matrices $A$ and $B$.  If $p(t) = \text{\rm{Tr}}
\left[(A+tB)^m\right]$ has nonnegative coefficients whenever 
$A,B$ are singular $n \times n$ 
positive semidefinite matrices, then $p(t)$ has nonnegative
coefficients whenever $A$ and $B$ are arbitrary $n \times n$
positive semidefinite matrices.
\end{thm}

The organization of this paper is as follows.  In Section \ref{prelim}, we
recall some facts about Hurwitz products, and in Section \ref{varprelim}
we derive the two equations found in Theorem \ref{mainvarthmAB}.  Finally,
in Section \ref{mainthmproof}, we use these equations
to prove our main Theorems \ref{maxtraceconj}, \ref{mainthm}, 
\ref{poseigvaluethm}, and \ref{singthm}.

\section{Preliminaries}\label{prelim}

We begin with a review of some basic facts involving Hurwitz products; some
of this material can be found in  \cite{JH3}.
The coefficients $S_{m,k}(A,B)$ may be generated via the
recurrence: 
\begin{equation}\label{basicrecurrence}
\begin{split}
S_{m,k}(A,B) = & \ AS_{m-1,k}(A,B)  + BS_{m-1,k-1}(A,B). \\
\end{split}
\end{equation}
The following lemmas will be useful for
computing the traces of the $S_{m,k}$.  
\begin{lem}\label{easytracelemma}
Fix integers $m > k \geq 0$.  For any two $n \times n$ matrices $A$ and $B$, we have
\[\text{\rm{Tr}}\left[S_{m,k}(A,B)\right] = \frac{m}{m-k}
\text{\rm{Tr}}\left[AS_{m-1,k}(A,B)\right].\]
\end{lem}
\begin{proof}
Consider the following chain
of equalities:
\begin{equation*}
\begin{split}
0 = & \ \text{Tr}\left[ {\sum\limits_{i = 1}^m {\left( {A + tB}
\right)^{i - 1} \left( {A - A} \right)\left( {A + tB} \right)^{m -
i} } } \right] \hfill \\
= & \ \text{Tr}\left[ {m A\left({A + tB} \right)^{m - 1}}\right]  -
\text{Tr} \left[\sum\limits_{i = 1}^m {\left( {A + tB} \right)^{i
- 1} A\left( {A + tB} \right)^{m - i} }  \right] \hfill \\
= & \ \text{Tr}\left[ {mA\left( {A + tB} \right)^{m - 1} } \right] -
\left. {\text{Tr}\left[ {\frac{d} {{dy}}\left( {Ay + tB} \right)^m
} \right]\;} \right|_{y = 1}  \hfill \\
= & \ \text{Tr}\left[ {mA\left( {A + tB} \right)^{m - 1} } \right] -
\left. {\frac{d} {{dy}}\left[ {\text{Tr}\left( {Ay + tB} \right)^m
} \right]\;} \right|_{y = 1}. \\
\end{split}
\end{equation*}
Since $S_{m,k}(Ay,B) = y^{m-k}S_{m,k}(A,B)$, it follows that the
coefficient of $t^k$ in the last expression above is just
\[m\text{Tr}[AS_{m-1,k}(A,B)] - (m-k)\text{Tr}[S_{m,k}(A,B)].\]
This proves the lemma.
\end{proof}

\begin{lem}\label{easytracelemma2}
Fix integers $m \geq k >0$.  For any two $n \times n$ matrices $A$ and $B$, we have
\[\text{\rm{Tr}}\left[S_{m,k}(A,B)\right] = \frac{m}{k}
\text{\rm{Tr}}\left[BS_{m-1,k-1}(A,B)\right].\]
\end{lem}
\begin{proof}
Follows from Lemma \ref{easytracelemma} by
taking the trace of both sides of equation (\ref{basicrecurrence}).
\end{proof}

Let $A$ and $B$ be $n \times n$ 
Hermitian matrices.  Since $S_{m,k}(A,B)$ is the sum of all words of 
length $m$ in $A$ and $B$ with $k$ $B$'s, it follows that the conjugate 
transpose of $S_{m,k}(A,B)$ simply permutes its constituent summands.   
This verifies the following fact.

\begin{lem}\label{realcoefflem}
If $A$ and $B$ are $n \times n$ Hermitian matrices, then 
the matrix $S_{m,k}(A,B)$ is Hermitian.
\end{lem}

\begin{cor}\label{realcoeffcor}
The polynomial $p(t) = \text{\rm{Tr}}\left[(A+tB)^m\right]$ has all
real coefficients whenever $A$ and $B$ are $n \times n$ Hermitian
matrices.
\end{cor}

Although $S_{m,k}(A,B)$ is Hermitian for Hermitian
$A$ and $B$, it need not be positive definite even when $A$ and $B$ 
are $n \times n$ positive definite matrices, $n > 2$.
Examples are easily generated, and computational experiments suggest that it is
usually not positive definite.  
Finally, we record a useful fact about positive definite congruence.  

\begin{lem}\label{pdconglem}
Let $C$ be any complex $n \times n$ matrix and let $A$ be an $n \times n$
positive semidefinite matrix.  Then $CAC^*$ is positive semidefinite.
\end{lem}
\begin{proof}
See \cite[p. 399]{HJ1}.
\end{proof}

\section{Derivation of the Euler-Lagrange Equations}\label{varprelim}

The arguments for our main theorems are based on a
variational observation.  It says that an expression Tr[$S_{m,k}(A,B)]$
is minimized or maximized when $A$ and $B$ satisfy the Euler-Lagrange 
equations (see Corollary \ref{mainvarcorAB}).
Before presenting a proof of this fact, we give a series of technical preliminaries.

\begin{prop}\label{mainvarpropAcomplex}
Let $m > k > 0$ be positive integers.
Fix $B$ to be any Hermitian $n \times n$ matrix, and suppose that $A$ is a positive 
semidefinite matrix of norm $1$ that minimizes (resp. maximizes) 
\[ \text{\rm Tr}[S_{m,k}(A,B)]\]
over all positive semidefinite matrices of norm $1$.
Let $\varepsilon > 0$, and let $C := C(x) = (c_{rs}(x))$ 
be an $n \times n$ matrix with entries $c_{rs}(x) = u_{rs}(x)+iv_{rs}(x)$ in which $u_{rs}$ and $v_{rs}$
are differentiable functions $u_{rs},v_{rs}: [-\varepsilon,\varepsilon] \to \mathbb R$.
Moreover, suppose that $C(0) = I$ and $CAC^* \neq 0$ for all 
$x \in [-\varepsilon,\varepsilon]$.  Then, the following  
identity holds: \[ \left. {\text{\rm Tr}\left[\frac{d}{dx}\left(\frac{CAC^*}{\|CAC^*\|} \right) 
S_{m-1,k}(A,B)\right] } \right|_{x = 0} = 0.\]
\end{prop}

\begin{proof}
Let $A, B,$ and $C$ be as in the statement of the theorem.  
Keeping in mind Corollary \ref{realcoeffcor}, we may
consider the differentiable function $f: [-\varepsilon,\varepsilon] \to
\mathbb R$ given by
\[f(x) = \text{Tr}\left[S_{m,k}\left(\frac{CAC^*}{\|CAC^*\|},B \right)\right].\]
By hypothesis, the minimum (resp. maximum) of $f$ is achieved at $x = 0$.  
Consequently, it follows that
\begin{equation}\label{derivzero}
\left. {\frac{{df(x)}} {{dx}}\;} \right|_{x = 0}  = 0.
\end{equation}
Next, notice that, \[\frac{d} {{dx}} {\text{Tr}\left[ \left(
{\frac{CAC^*}{\|CAC^*\|} + tB} \right)^m \right] } =
\text{Tr}\left[ {\frac{d} {{dx}}\left( {\frac{CAC^*}{\|CAC^*\|}
+ tB} \right)^m } \right]
\]
\[=  \text{Tr}\left[ {\sum\limits_{i = 1}^m {\left( {\frac{CAC^*}{\|CAC^*\|}
 + tB} \right)^{i - 1} \frac{d} {{dx}}\left( {\frac{CAC^*}{\|CAC^*\|}
 + tB} \right)\left( {\frac{CAC^*}{\|CAC^*\|} +
tB} \right)^{m - i} } } \right].\]  
When $x= 0$, the above expression evaluates to
\[\left. { \text{Tr}\left[ {\sum\limits_{i = 1}^m {\left( {A+ tB} \right)^{i - 1} 
\frac{d} {{dx}}\left( {\frac{CAC^*}{\|CAC^*\|}} \right)\left( {A+
tB} \right)^{m - i} } } \right]}  \right|_{x = 0} \]
\[= \left. {\text{Tr}\left[ m \frac{d}{dx} \left( \frac{CAC^*}{\|CAC^*\|} \right) \left( A + tB \right)^{m - 1}\right]
}  \right|_{x = 0} \]
It follows, therefore, from (\ref{derivzero}) that
\[ \left. {\text{\rm Tr}\left[\frac{d}{dx}\left(\frac{CAC^*}{\|CAC^*\|} \right) 
S_{m-1,k}(A,B)\right] } \right|_{x = 0} = 0.\]
This completes the proof.
\end{proof}

A corresponding statement can be made by fixing $A$ and 
minimizing (resp. maximizing) over $B$.

\begin{prop}\label{mainvarpropBcomplex}
Let $m > k > 0$ be positive integers.
Fix $A$ to be any Hermitian $n \times n$ matrix, and let $B$ be a positive 
semidefinite matrix of norm $1$ that minimizes  (resp. maximizes) \[ \text{\rm Tr}[S_{m,k}(A,B)]\]
over all positive semidefinite matrices of norm $1$.
Let $\varepsilon > 0$, and let $C := C(x) = (c_{rs}(x))$ 
be an $n \times n$ matrix with entries $c_{rs}(x) = u_{rs}(x)+iv_{rs}(x)$ in which $u_{rs}$ and $v_{rs}$
are differentiable functions $u_{rs},v_{rs}: [-\varepsilon,\varepsilon] \to \mathbb R$.
Moreover, suppose that $C(0) = I$ and $CBC^* \neq 0$ for all 
$x \in [-\varepsilon,\varepsilon]$.  Then, the following  
identity holds: \[ \left. {\text{\rm Tr}\left[\frac{d}{dx}\left(\frac{CBC^*}{\|CBC^*\|} \right) 
S_{m-1,k-1}(A,B)\right] } \right|_{x = 0} = 0.\]
\end{prop}

\begin{proof}
The proof is similar to that of Proposition \ref{mainvarpropAcomplex}, so 
we omit it.
\end{proof}

In our next lemma, we compute the derivative 
found in Propositions \ref{mainvarpropAcomplex} and 
\ref{mainvarpropBcomplex}.  For notational simplicity, the entry-wise 
derivative of the matrix $C$ evaluated at the point $x = 0$ 
will be denoted by $C'$.

\begin{lem}\label{derivevallemma}
With the hypotheses as in Proposition \ref{mainvarpropAcomplex},
we have
\[ \left. {\frac{d}{dx} \left(\frac{CAC^*}{\left\|CAC^*\right\|} \right)} \right|_{x = 0}
=  C'A+ AC'^*   -  \text{\rm Tr}[C'A^2]A -  \overline{\text{\rm Tr}[C'A^2]}A. \]
\end{lem}
\begin{proof}
A straightforward application of the product rule \cite{HJ2}
for (matrix) differentiation gives us that 
\begin{equation}
\begin{split}
\frac{d}{dx}\left(\frac{CAC^*}{\left\|CAC^*\right\|} \right) = \ 
& \frac{d}{dx}\left(\frac{1}{\left\|CAC^*\right\|} \right)CAC^* +
\frac{1}{\left\|CAC^*\right\|}\left(\frac{dC}{dx}AC^*  + C A\frac{dC^*}{dx}\right). \\
\end{split}
\end{equation}
Next, we compute that
\begin{equation*}
\begin{split}
\frac{d}{dx} \left\|CAC^*\right\|^{-1}   = & \ 
- \left\|CAC^*\right\|^{-2}\frac{d}{dx} \left\|CAC^*\right\|  \\
 = & \ - \left\|CAC^*\right\|^{-2}\frac{d}{dx}\left( \text{Tr}[CAC^*CAC^*]^{1/2}  \right)  \\
  = & \ -(1/2) \left\|CAC^*\right\|^{-2}\text{Tr}[CAC^*CAC^*]^{-1/2} \frac{d}{dx} \text{Tr}[CAC^*CAC^*]. \\
\end{split}
\end{equation*}

The product expansion of $\frac{d}{dx} \text{Tr}[CAC^*CAC^*]$ occurring 
in this last line is:
\begin{equation*}
\begin{split}
& \ \text{Tr}\left[\frac{dC}{dx}AC^*CAC^*+CA\frac{dC^*}{dx}CAC^*+
CAC^*\frac{dC}{dx}AC^*+CAC^*CA\frac{dC^*}{dx}\right]   \\
= & \ 2\text{Tr}\left[\frac{dC}{dx}AC^*CAC^*\right]+2\overline{\text{Tr}\left[\frac{dC}{dx}AC^*CAC^*\right]}.
\end{split}
\end{equation*}
Finally, setting $x = 0$ and using the assumptions that $\|A\| = 1$
and $C(0) = I$, it follows that 
\[  \left. { \frac{d}{dx}\left(\frac{CAC^*}{\left\|CAC^*\right\|} \right) } \right|_{x = 0}
 =  C'A + AC'^*- \text{Tr}[C'A^2]A -  \overline{\text{Tr}[C'A^2]}A,\]
completing the proof of the lemma.
\end{proof}

We now have enough to prove the main results of this section.

\begin{thm}\label{mainvarthmcomplexA}
Let $m > k > 0$ be positive integers.
Fix $B$ to be any Hermitian $n \times n$ matrix, and let $A$ be a positive 
semidefinite matrix of norm $1$ that minimizes  (resp. maximizes) 
\[ \text{\rm Tr}[S_{m,k}(A,B)]\]  
over all positive semidefinite matrices of norm $1$.
Then, the following identity holds:  
\[AS_{m-1,k}(A,B)= A^2 \text{\rm Tr}[AS_{m-1,k}(A,B)].\]
\end{thm}
\begin{proof}
Let $A$ and $B$ be as in the hypotheses of the theorem.
By using different matrices $C$ in the statement of Proposition 
\ref{mainvarpropAcomplex}, we will produce a set of equations satisfied by the 
entries of $AS_{m-1,k}(A,B)$ that combine to 
make the single matrix equation above.  
For ease of presentation, we introduce the following
notation.  For integers $r,s$, let $E_{rs}$ denote the $n \times n$ matrix with
all zero entries except for a $1$ in the $(r,s)$ entry.

Fix integers $1 \leq r,s \leq n$ and take $C = I + xE_{rs}$.  Since $C$ is invertible for
all $x \in [-1/2,1/2]$, it follows that $CAC^* \neq 0$ for all such $x$.  
Therefore, the hypotheses of
Lemma \ref{derivevallemma} are satisfied.  The formula 
there gives us that
\begin{equation}\label{deriveval}
\begin{split}
\left. {\frac{d}{dx}\left(\frac{CAC^*}{\|CAC^*\|} \right)}  \right|_{x = 0}
 = & \ C'A+ AC'^*   -  \text{\rm Tr}[C'A^2]A -  \overline{\text{\rm Tr}[C'A^2]}A. \\
\end{split}
\end{equation}
Additionally, Proposition \ref{mainvarpropAcomplex}, along with
a trace manipulation, tells us that 
\begin{equation}\label{mainvarpropconseq}
\begin{split}
 \ \left(\text{Tr}[C'A^2]+ \overline{\text{Tr}[C'A^2]} \right) &  \text{\rm Tr}\left[AS_{m-1,k}(A,B) \right] \\
= & \ \text{\rm Tr}\left[C'AS_{m-1,k}(A,B) + AC'^*S_{m-1,k}(A,B) \right]  \\
= & \ 
 \text{\rm Tr}\left[C'AS_{m-1,k}(A,B) \right]
+\text{\rm Tr}\left[S_{m-1,k}(A,B)AC'^* \right] \\
= & \  \text{\rm Tr}\left[C'AS_{m-1,k}(A,B) \right]
+\overline{\text{\rm Tr}\left[C'AS_{m-1,k}(A,B) \right]}.
\end{split}
\end{equation}
Since $C' = E_{rs}$, a computation shows that for any matrix $N$,
the trace of $C'N$ is just the $(s,r)$ entry of $N$.  In particular, it follows from
(\ref{mainvarpropconseq}) that the $(s,r)$ entries of 
$AS_{m-1,k}(A,B)+ \overline{AS_{m-1,k}(A,B)}$ 
and $(A^2+ \overline{A^2}) \text{\rm Tr}[AS_{m-1,k}(A,B)]$ coincide.
We have therefore proved the following identity of matrices:
\begin{equation}\label{realeqidentity}
AS_{m-1,k}(A,B)+ \overline{AS_{m-1,k}(A,B)}= (A^2+ \overline{A^2})\text{\rm Tr}[AS_{m-1,k}(A,B)].
\end{equation}

We next perform a similar examination using the matrices $C = I+ixE_{rs}$
to arrive at a second matrix identity.  
Combining equation (\ref{deriveval}) and Proposition \ref{mainvarpropAcomplex}
as before, we find that 
\begin{equation}\label{imaginaryidentity}
AS_{m-1,k}(A,B) - \overline{AS_{m-1,k}(A,B)}= (A^2 - 
\overline{A^2} ) \text{\rm Tr}[AS_{m-1,k}(A,B)].
\end{equation}
The theorem now follows by adding these two equations and 
dividing both sides of the result by $2$.
\end{proof}

Similar arguments using Proposition \ref{mainvarpropBcomplex} in place
of Proposition \ref{mainvarpropAcomplex} produce the following results.

\begin{thm}\label{mainvarthmcomplexB}
Let $m > k > 0$ be positive integers.
Fix $A$ to be any Hermitian $n \times n$ matrix, and let $B$ be a positive 
semidefinite matrix of norm $1$ that minimizes  (resp. maximizes) 
\[ \text{\rm Tr}[S_{m,k}(A,B)]\]
over all positive semidefinite matrices of norm $1$.
Then, the following identity holds:  
\[BS_{m-1,k-1}(A,B)= B^2 \text{\rm Tr}[AS_{m-1,k-1}(A,B)].\]
\end{thm}

Combining the statements of this section, we have finally derived 
the Euler-Lagrange equations (\ref{matrixeqs}) for
Conjecture \ref{BMV}.

\begin{cor}\label{mainvarcorAB}
Let $m > k > 0$ be positive integers, 
and let $A$ and $B$ be positive semidefinite
matrices of norm $1$ that minimize  (resp. maximize) 
the quantity $\text{\rm Tr}[S_{m,k}(A,B)]$
over all positive semidefinite matrices of norm $1$.
Then $A$ and $B$ satisfy the following pair of equations:
\begin{equation*}
\begin{cases}
AS_{m-1,k}(A,B) &= \ A^2 \text{\rm Tr}[AS_{m-1,k}(A,B)]  \\
BS_{m-1,k-1}(A,B) &= \ B^2 \text{\rm Tr}[BS_{m-1,k-1}(A,B)].
\end{cases}
\end{equation*}
\end{cor}

We remark that our proof generalizes directly to show 
that the same equations hold for Hermitian minimizers 
(resp. maximizers), or more generally, for classes of unit norm matrices
with the same inertia.  This result is the main ingredient in our proof
of Theorem \ref{maxtraceconj} concerning the maximum of $\text{Tr}[S_{m,k}(A,B)]$.

\begin{cor}\label{hermitianeqs}
Let $m > k > 0$ be positive integers, and
let $A$ and $B$ be Hermitian matrices of 
norm $1$ that minimize  (resp. maximize) the quantity $\text{\rm Tr}[S_{m,k}(A,B)]$
over all Hermitian matrices of norm $1$.
Then $A$ and $B$ must satisfy the following pair of equations:
\begin{equation*} 
\begin{cases}
AS_{m-1,k}(A,B) & = \ A^2 \text{\rm Tr}[AS_{m-1,k}(A,B)]  \\
BS_{m-1,k-1}(A,B) & = \ B^2 \text{\rm Tr}[BS_{m-1,k-1}(A,B)].
\end{cases}
\end{equation*}
\end{cor}

In general, we conjecture that trace minimizers commute (Conjecture \ref{commuteconj}), a
claim that would imply Conjecture \ref{BMV}.  
We close this section with one more application of the Euler-Lagrange equations.

\begin{cor}
Suppose that the minimum of $\text{\rm Tr}[S_{m,k}(A,B)]$ over the
set of positive semidefinite matrices is 
zero, and let $A$ and $B$ be positive semidefinite matrices
achieving this minimum.  Then, $S_{m,k}(A,B) = 0$.
\end{cor}

\begin{proof}
When $k = m$ or $k = 0$, the claim is clear.  Therefore,
suppose that $m > k > 0$.
Let $A$ and $B$ be positive semidefinite matrices with
$\text{\rm Tr}[S_{m,k}(A,B)] = 0$.  If either of $A$ or $B$ is zero
then the corollary is trivial.  Otherwise, consider
\[ 0 = \frac{\text{\rm Tr}[S_{m,k}(A,B)]  }{\|A\|^{m-k}\|B\|^{k}}= 
\text{\rm Tr}[S_{m,k}(\widetilde{A},\widetilde{B})],\]
in which $\widetilde{A} = A/\|A\|$ and  $\widetilde{B} = B/\|B\|$.
Combining equations (\ref{matrixeqs})
with the assumptions, it follows that $\widetilde{A}S_{m-1,k}(\widetilde{A},\widetilde{B}) = 0$ 
and $\widetilde{B}S_{m-1,k-1}(\widetilde{A},\widetilde{B}) = 0$.  
Moreover, equation (\ref{basicrecurrence}) implies that
\[S_{m,k}(\widetilde{A},\widetilde{B}) = \widetilde{A}S_{m-1,k}(\widetilde{A},\widetilde{B}) +\widetilde{B}S_{m-1,k-1}(\widetilde{A},\widetilde{B}) = 0.\]
Multiplying both sides of this identity by $\|A\|^{m-k}\|B\|^{k}$ 
completes the proof.
\end{proof}

\section{Proofs of the Main Theorems}\label{mainthmproof}

We first use the Euler-Lagrange equations to prove Theorem \ref{maxtraceconj}.

\begin{proof}[Proof of Theorem \ref{maxtraceconj}]
Let $m > 1$ and $n$ be positive integers.
Since our arguments are the same in both cases,
we consider determining the maximum
of $\text{\rm Tr}[S_{m,k}(A,B)]$. 
Let $M$ be the compact set of Hermitian
matrices with norm 1 and choose $(A,B) \in M\times M$
that maximizes $\text{Tr}[S_{m,k}(A,B)]$. 
If $k = 0$, then the desired inequality
is of the form \[ \text{Tr}[A^m] \leq \sum_{i=1}^n |\lambda_i|^m 
\leq \sum_{i=1}^n \lambda_i^2 = \|A\| = 1,\]
in which $\lambda_1,\ldots,\lambda_n$ are the
eigenvalues of $A$.  A similar argument holds
for $m= k$.  Therefore, we assume below that $m > k > 0$.

The Euler-Lagrange equations from Corollary \ref{hermitianeqs}
imply that
\begin{equation}\label{hermeq}
 AS_{m-1,k}(A,B) = A^2 \text{Tr}[AS_{m-1}(A,B)].
\end{equation}
Performing a uniform, unitary similarity, we may 
assume that $A$ is diagonal of the form 
$A = \text{diag}(\lambda_1,\ldots,\lambda_r,0,\ldots,0)$, in which
$\lambda_1,\ldots,\lambda_r$ are nonzero.  
Let $\widetilde{A} = \text{diag}(\lambda_1^{-1},\ldots,\lambda_r^{-1},0,\ldots,0)$ 
be the pseudo-inverse of $A$, and set $D = \widetilde{A}A$. 
Multiplying both sides of \ref{hermeq} by $\widetilde{A}$, it follows that 
\[ DS_{m-1,k}(A,B) = A\text{Tr}[AS_{m-1,k}(A,B)].\]
Taking the norm of both sides of this expression and applying
Lemma \ref{easytracelemma}, we have
\begin{equation}\label{hermproofchain}
 \frac{m-k}{m}\text{Tr}[S_{m,k}(A,B)] = 
\|DS_{m-1,k}(A,B)\| \leq \|S_{m-1,k}(A,B)\| \leq {m-1 \choose k}.
\end{equation}
It follows that ${m \choose k} \leq \text{Tr}[S_{m,k}(A,B)]  \leq 
\frac{m}{m-k} {m-1 \choose k} = {m \choose k}$ as desired.

We next verify the final assertions in the 
statement of the theorem.
From above, every
inequality in the chain (\ref{hermproofchain}) is an actual
equality.  Thus, each term occurring in 
\[\|S_{m-1,k}(A,B)\| = \sum_{W} {\|W(A,B)\|},\]
a sum over length $m-1$ words $W$ with $k$ $B$'s,
takes the value $1$.  In particular, we have that 
$1 = \|A^{m-k-1}B^k\|$.  When $m-1 > k > 1$, an application of Lemma
\ref{normeqlem} below completes the proof of the theorem.
The remaining cases $k = 1$ or $m = k+1$ are dealt with as follows.

Without loss of generality, we may suppose that $k = 1$ (interchange
the roles of the matrices $A$ and $B$).
Applying the Cauchy-Schwartz inequality, we obtain the following chain of inequalities,
\[ 1 = \text{Tr}[A^{m-1}B]^2 = \left(\sum_{i=1}^n {\lambda_i^{m-1}b_{ii}}\right)^2
\leq \sum_{i=1}^n {(\lambda_i^2)^{m-1}}\sum_{i=1}^n {b_{ii}^2} \leq  \|A\| \|B\| = 1.\]
It follows that each inequality above is actually an equality.  In particular,
the second to last identity says that $B$ is diagonal.  Moreover,
equality in Cauchy-Schwartz implies that $\lambda_i^{m-1} = \delta b_{ii}$
for some real number $\delta$ and all $i$.  Since $1 = \left|\sum_{i=1}^n \lambda_i^{m-1}b_{ii}\right| = \delta^2$,
it follows that $A = \pm B$.
If in addition, $m > 2$, and $A$ has more than $1$ nonzero
eigenvalue, then  \[1 = \sum_{i=1}^n \lambda_i^2 > 
\sum_{i=1}^n (\lambda_i^2)^{m-1} = \sum_{i=1}^n b_{ii}^2 = 1,\] a contradiction.
Therefore, the conclusions of the theorem hold  for $k = 1$, and this
finishes the proof.
\end{proof}

\begin{lem}\label{normeqlem}
Suppose that $A$ and $B$ are Hermitian matrices of norm $1$
and $r > 0$ and $s > 1$ are integers such that $\|A^{r}B^s\| = 1$.
Then, $A = \pm B$ has only $1$ nonzero eigenvalue.
\end{lem}
\begin{proof}
Performing a uniform, unitary similarity, we may suppose that
$B$ is a diagonal matrix  with entries less than or equal to $1$
in absolute value.
From the hypotheses, we have that 
\[1 = \|A^{r}B^s\| \leq \|A^{r}\|\|B^s\| \leq  \|A\|^{r}\|B\|^s = 1.\]
Therefore, $\|B^s\| = 1 = \|B\|$, and since $s > 1$, this implies
that $B$ has a single nonzero eigenvalue.  
It follows that $\|A^{r}B^s\| = \|A^{r}B\| = 1$ is equal to the absolute
value of the $(1,1)$ entry of $A^r$.  Finally, since $\|A^r\| = 1$, the
matrix $A^r$ has only one nonzero entry, and therefore, $A$ has
only one nonzero eigenvalue.  
Thus, $A^r = \pm A$ and since $A^r = \pm B$, it follows that $A = \pm B$.
This completes the proof.
\end{proof}

The argument for our next result 
requires the following well-known fact; we 
provide a proof for completeness.

\begin{lem}\label{prodpdlemma}
If $P$ and $Q$ are positive semidefinite matrices,
then $PQ$ has all nonnegative eigenvalues.
\end{lem}
\begin{proof}
Suppose first that $P$ is positive definite.  Then $PQ$ is
similar to \[P^{-1/2}PQP^{1/2} = P^{1/2}QP^{1/2}.\]
In particular, $PQ$ is similar to a positive semidefinite matrix
by Lemma \ref{pdconglem}.  Therefore, in this case $PQ$ 
has all nonnegative eigenvalues. 
The general version of the claim now follows from continuity.
\end{proof}

We are now prepared to present a proof
that counterexamples to Conjecture \ref{BMV}
are closed upward.  Theorem \ref{poseigvaluethm} 
closely follows.

\begin{proof}[Proof of Theorem \ref{mainthm}]
Suppose that Conjecture \ref{BMV} is false for some
$m$ and $k$ and let $A$ and $B$ be real positive semidefinite matrices
of unit norm that minimize
\[\text{Tr}[S_{m,k}(A,B)] = \frac{m}{m-k}\text{Tr}[AS_{m-1,k}(A,B)] = 
\frac{m}{k}\text{Tr}[BS_{m-1,k-1}(A,B)].\] 
We show that for these same matrices $A$ and $B$, 
we have $\text{Tr}[S_{m+1,k}(A,B)] < 0$ and $\text{Tr}[S_{m+1,k+1}(A,B)] < 0$.  

Combining equation (\ref{basicrecurrence}) and 
using the identities (\ref{matrixeqs})
from Corollary \ref{mainvarcorAB}, it follows that 
\begin{equation}\label{negdefeq}
S_{m,k}(A,B) = A^2 \text{Tr}[AS_{m-1,k}(A,B)] + B^2 \text{Tr}[BS_{m-1,k-1}(A,B)].
\end{equation}
This matrix is negative semidefinite since it is the
sum of two such matrices.  Hence, the product $AS_{m,k}(A,B)$ has all 
non-positive eigenvalues by Lemma \ref{prodpdlemma}.   Thus, Lemma
\ref{easytracelemma} implies that 
\begin{equation}\label{inductivetreq}
\text{Tr}[S_{m+1,k}(A,B)] = \frac{m}{m+1-k} \text{Tr}[AS_{m,k}(A,B)] \leq 0.
\end{equation}
In the case of equality, 
multiplying equation (\ref{negdefeq}) on the left by $A$ and 
taking the trace of both sides, it follows that  
\[0 = \text{Tr}[A^3] \text{Tr}[AS_{m-1,k}(A,B)] + \text{Tr}[AB^2] \text{Tr}[BS_{m-1,k-1}(A,B)].\]
However, Tr$[AB^2] \geq 0$ by Lemma \ref{prodpdlemma} and 
since $A$ is nonzero, we must have $\text{Tr}[A^3] > 0$.  This
gives a contradiction to equality in (\ref{inductivetreq}).  
It follows that \mbox{$\text{Tr}[S_{m+1,k}(A,B)]  < 0$} as desired.  

In the same manner, we can also prove that 
\[\text{Tr}[S_{m+1,k+1}(A,B)] = \frac{m+1}{k+1} \text{Tr}[BS_{m,k}(A,B)]\]
is negative.  The conclusions of the theorem now follow immediately.
\end{proof}

\begin{proof}[Proof of Theorem \ref{poseigvaluethm}]
Suppose that \text{\rm Tr}$[S_{m,k}(A,B)]$ can be made negative.  
The proof of Theorem \ref{mainthm} shows that there exist
positive semidefinite matrices $A$ and $B$ such that
$S_{m,k}(A,B)$ is negative semidefinite and
 \mbox{Tr$[S_{m,k}(A,B)] < 0$} (so that $S_{m,k}(A,B)$ is nonzero).
It follows that the second implication in the statement of the theorem is false.
The converse is clear.
\end{proof}

Finally, we work out the proof of Theorem \ref{singthm}; the
argument is similar in spirit to the proof of Theorem \ref{maxtraceconj}.

\begin{proof}[Proof of Theorem \ref{singthm}]
Suppose we know that Conjecture \ref{BMV} is true for
the power $m-1$ and also suppose that for some $k$
there exist $n \times n$ positive definite matrices $A$ and $B$
such that $\text{Tr}[S_{m,k}(A,B)]$ is negative.  
Clearly, we must have $m > k > 0$.
By homogeneity, there are positive definite $A$
and $B$ with norm 1 such that $\text{Tr}[S_{m,k}(A,B)]$ is
negative.  Let $M$ be the set of positive
semidefinite matrices with norm 1 and choose $(A,B) \in M\times M$
that minimizes $\text{Tr}[S_{m,k}(A,B)]$;  our goal is
to show that $A$ and $B$ must both be singular.

Suppose by way of contradiction that $A$ is invertible.  
The Euler-Lagrange equations say that
\[ AS_{m-1,k}(A,B) =  A^2 \text{\rm Tr}[AS_{m-1,k}(A,B)].\]
Multiplying both sides of this equation by $A^{-1}$
and taking the trace, it follows
that \[\text{Tr}[S_{m-1,k}(A,B)] = \text{Tr}[A] \text{\rm Tr}[AS_{m-1,k}(A,B)].\]
By hypothesis, $\text{Tr}[S_{m-1,k}(A,B)]$ is nonnegative.  Therefore,
using Lemma \ref{easytracelemma}, we have
\begin{equation*}
\begin{split}
\frac{m-k}{m}\text{\rm{Tr}}\left[S_{m,k}(A,B)\right] =
& \ \text{\rm{Tr}}\left[AS_{m-1,k}(A,B)\right] \\
= & \  \frac{\text{Tr}[S_{m-1,k}(A,B)]}{\text{Tr}[A]} \\
 \geq & \ 0,
\end{split}
\end{equation*}
a contradiction (Tr$[A]$ is nonzero since $A$ is nonzero).  
It follows that $A$ must be singular.
A similar examination with $B$ also shows that it must
be singular.

Thus, if Conjecture \ref{BMV} is true for singular $A$ and $B$,
it must be true for invertible $A$ and $B$ as well.  This completes
the proof of the theorem.
\end{proof}

\section{Acknowledgments}

We would like to thank Scott Armstrong for several interesting and useful
discussions concerning a preliminary version of this manuscript.


\end{document}